
\magnification 1200
\input amstex
\documentstyle{amsppt}
\def\ju{\vskip .8truecm plus .1truecm minus .1truecm}

\def\sju{\vskip .4truecm plus .1truecm minus .1truecm}
\def\lra{\longrightarrow}
\def\lqq{\lq\lq}
\def\se{\subseteq}
\def\shimply{\Rightarrow} 
\def\RR{\Bbb R}
\def\al{\alpha}
\def\ep{\epsilon}
\def\om{\omega}

\nologo
\NoBlackBoxes
\hsize=13.5truecm
\vsize=22.5truecm

\vskip 1 truecm

\centerline{\bf IS THE FUNCTION $\bold 1/\bold x$ CONTINUOUS AT 0?\footnote[]{This work has originated around year 2000. The author gave the first public lecture on the topic no later than February 7, 2003, at Texas A\&M University. Current paper is a minor revision of the version written on December 8, 2003. The paper has been in circulation since that date. \copyright\,\, R. Dimitri\'c\,\,\, 2003 }}

\vskip 0.5truecm
\centerline{by {\it Radoslav M. Dimitri\'c}} 
\vskip 0.3truecm
\centerline{email: {\it dimitricr\@member.ams.org}}

\vskip 1truecm
 
 {\bf Question Number Two}: Is the function 
$$   f(x)=\cases
0, & \text{ if } x=0\cr
1, & \text{ if } x\geq 1
\endcases    \hskip  4truecm  (1)$$   
continuous at $0$? 
If we open any elementary calculus textbook in use, 
the answer will almost invariably be \lqq no" to both
questions. Moreover with some textbooks, these are the only 
\lqq discontinuities" various functions seem to have.  
The daily calculus teaching practice follows the same party line. 
We will give correct answers to these questions in the subsequent pages. 

First let us briefly look into the development of the ideas and the concepts.
\sju
{\bf What is a continuous function?}

It is worth noting that the concept of a function as we know it today was
arrived at only through numerous historical meanderings. Historically, this concept has
been closely intertwined with that of continuity. [Ferraro, 2000], for instance,
gives a detailed account of the development of this concept with Euler.
A function for Euler was by definition continuous and differentiable and 
expandable into a Taylor series. The functions were intrinsically continuous,
for \lqq the variables varied in a continuous way". In the eighteenth century a
function was restricted to a quantity definable by a single formula, and Euler
followed the same practice. The following was not
considered to be a function:
$$  {
f(x)=
\cases
x^2, & \text{ if } x\leq 0,\cr
5x, & \text{ if } x> 0
\endcases     \hskip  3.9truecm  (2) 
}$$
By this definition, $1/x$ is a function, and a continuous one since it is defined
by a unique formula; the fact that its graph consisted of two pieces was of 
no consequence to Euler. By the same token the integer part function [\,\,] would
likely be considered to be  a continuous function by Euler. [Cauchy, 1844]
quickly points out some difficulties with Euler's view by noting that the
function $\sqrt{x^2}$ would be considered a continuous function by Euler,
but the same function written in two pieces: $=x$, if $x\geq 0$ and $=-x$, if $x\leq 0$
would not be considered to be a function. Only in 1837 did Dirichlet get rid
of this notion of function interpreted as a \lqq unique" formula.

The notion of  function has nowadays crystallized into a well-founded concept. 
A function
is no longer only a formula (or a compact formula at that), 
but is a triple $f:A\lra B$ where $A$ is the domain
for its variables, $f$ is a relation with unique second component value, 
for every
independent variable, and $B$ is the codomain (containing the range); see
here [Bourbaki, 1939], for example.
Thus, the functions should be taught properly: domain, codomain and the assignment
rule(s) with certain properties make up a function.
Precisely speaking, a function $f$ from $A$ to $B$ is the ordered triple
$(f,A,B)$, (also denoted by $f:A\lra B$), where $f\se A\times B$ is a relation
(i.e. a subset),
with the property that, for every $a\in A$ there is a unique $b\in B$ such that
$(a,b)\in f$. Thus, two functions are equal, if all three of their corresponding
parts (the formula, domain and codomain) are equal. One might say that using relations to define
functions may be a bit overreaching in a basic calculus course and perhaps it is;
on the other hand, I see textbooks on \lqq precalculus" discuss the subject of
relations. I encounter students' puzzlement, or even confusion when they hear
the \lqq new word" `codomain' (some also use the word `target' for the same
thing). They've heard of the range, but what is this?
Codomain is as important as the domain and is perfectly dual to it (and that is
a very important point). Given a formula, we at first simply do not know what the range
will be (in any but trivial cases), but can often give a \lqq rough" target. 

Consider the following example:
For a square black-and-white photograph, for every point on it (and one tends
to modernize terminology, not quite accurately and say \lqq to every pixel"), 
assign a value to it in the interval $[0,1]$
that denotes the amount of black at that point; thus 1 would be assigned to a perfectly black
point and 0 to a perfectly white. The interval $[0,1]$ is a natural codomain,
although we could go up to $\RR$ if we are only talking real functions; anything
in between can also be taken to be a codomain; as we change the codomains we
get different function each time, although the assignment rule (the formula) would stay
 the same. On the
other hand $[-1,0)$ would not be a valid codomain. We don't know the range
and it would be rather difficult to find -- very dark photos would have the range
in the upper half of the interval $[0,1]$, the light ones in the bottom half.
We know however that the range is the smallest of all the codomains. 
Dually, given a formula, we speak of the domain to fit the formula, whereas
we really mean the largest domain for which the formula works. Any subset of
the maximal domain can also be a domain. As we change these (without changing
the formula) we do get different functions. Thus, for the same formula
$f(x)=x^2$, we have various possibilities 

\noindent a) $f:\RR\lra\RR$, b) $f:\RR^+\lra \RR$, c) $f:\RR^+\lra\RR^+$, 
d)  $f:\RR\lra \RR^+$ \hfill (3)

\noindent which give us  different
functions with different properties: the first is neither one-to-one nor
onto (i.e. range is not equal to its codomain), the second one is one-to-one
but not onto while the function  in c) has an inverse (it is one-to-one and onto),
unlike the rest of them; the function in d) is onto but not one-to-one.
 In this way, not only is 
$1/x$ not a function at $0$, but 
we can point out to our students some crucial facts, such as that the property
of being one-to-one depends on the domain (such as in functions (3) above).
This then helps better understand the notion of
inverse function, the relationship between the range and the codomain, etc.\footnote{To be fair, 
better, or more advanced books on analysis often have a correct treatment of the notion of function,
and in particular that of an \lqq onto" function and a codomain; see here for instance [Rudin, 1976] or
[Marde\v si\'c, 1974]}

Now a definition of continuity:
Given a function $f:A\lra\RR$ on a real domain $A\se \RR$, it is 
{\it continuous} at $a\in A$, if\footnote{Here we do not require that the
limit is taken over a non-isolated point -- $x$'s may equal $a$,} 
$$   \lim_{x\to a}f(x)=f(\lim_{x\to a}x)=f(a).   \hskip  3.2truecm   (4) $$
Implicitly, this definition is in the form of an implication: If $a\in A$, then 
the function $f$ and the  limit as $x\to a$ can be commuted. This is fairly close
to how Bolzano defines
continuity in his privately published manuscript  on what we today
call the {\it intermediate value theorem} [Bolzano, 1817]:
\lqq A function $f(x)$ varies according to the law of continuity for all
values of $x$ which lie inside or outside certain limits, is nothing other
than this: If $x$ is any such value, the difference $f(x+\om)-f(x)$ can be made 
smaller than any given quantity, if one makes $\om$ as small as one wishes."
For [Euler, 1748], this is a consequence of what he considers to be a function. 

Years later, [Weierstrass, 1874] gave a similar, somewhat more formal definition:
\lqq Here we call a quantity $y$ a continuous function of $x$, if upon taking a
quantity $\ep$, the existence of $\delta$ can be proved, such that for any
value between $x_0-\delta...x_0+\delta$, the corresponding value of $y$ lies
between $y_0-\ep...y_0+\ep$."  
(The difference
$f(x)-f(x_0)$ can be made arbitrarily small, if the difference $x-x_0$ is made
sufficiently small.)

In our modern quantifier notation, $f$ is continuous at $x_0\in \text{Dom} f$, if
$$   { \forall\ep>0\exists\delta>0\forall x\in A (|x-x_0|<\delta\shimply 
|f(x)-f(x_0)|<\ep).      \hskip 0.5truecm   (5)}$$ 

It is not clear whether Bolzano's paper had been well known in his own time
and in particular whether Cauchy was familiar with the content of this paper; 
see [Grattan-Guinness, 1969] on this. Whatever the answer, [Cauchy, 1821, p.43]
gives his own definition of continuity:

\lqq ... $f(x)$ will be called a {\it continuous} function, if...the numerical
values of the difference $f(x+\al)-f(x)$ decrease indefinitely with those 
of $\al$..." (infinitesimally small changes in $x$ should 
lead to infinitesimally small changes in $f$).

On a slightly more fancy level, if $f:A\lra B$ 
is a function between two metric spaces, then it is continuous at $a\in A$,
if, for every (open) ball $V$ in $B$ centered at $f(a)$ there is an (open) ball
centered at $a$ completely mapped into $V$.

Thus, the  mid-eighteenth century notion of \lqq continuity" referred to 
uniformity (wholeness) 
of the formula used to define the function; piecewise defined functions, 
such as (2), 
were not deemed to be functions (or to be \lqq continuous") under this notion. 
D'Alembert was one of the
proponents of defining (continuous) functions in this restricted sense. 
This was challenged by
 Daniel Bernoulli, after d'Alembert's work in 1747 on the motion
of a vibrating string (d'Alembert, 1747;
 see here for instance [Struik, 1969, pp.351--368] 
and Truesdell's introduction to [Euler, 1960] titled: 
\lqq The rational mechanics of 
flexible or elastic bodies, 1638--1788"). 
D'Alembert's solution to the partial differential equation describing motion
of an elastic vibrating string was of the form: $z(x,t)=f(t+x)+F(t-x)$, which was
unusual in that the solution was a combination of two arbitrary functions, thus
could not be considered to be a (continuous) function. [Euler, 1765] resolved this problem by
change of terminology and interpretations of the constants that come out in 
integrating the given pde's [Feraro, 2000].

The alternative notion of \lqq continuity"
that was coming into greater prominence referred to functions that
\lqq can be produced by a free motion of the hand."
This notion of being able to draw the graph of a function without lifting 
a pencil is
that of {\it contiguity} as formulated in 1791 by Louis Arbogast: \lqq
The law of continuity is again broken when the different parts of a curve do
not join to one another... We will call curves of this kind discontiguous curves,
because all their parts are not contiguous, and similarly for discontiguous 
functions" (see [Jourdain, 1913], pp.675--676). The interest in 
(dis)continuities was further heightened by Fourier in his celebrated 
work on heat [Fourier, 1822]. A definite disturbance regarding the concept of
continuity and a need to fix the
concept was induced by Peter Lejeune-Dirichlet who gave, in 1829, his example
of the function equal to a constant over the set of rational numbers and
equal to another constant on the set of the irrationals; this function was
not continuous even at a single point.\footnote{Another nice example was given
by Riemann: $f:[0,1]\lra [0,1]$ is defined to be zero on irrational numbers and $=1/q$, for 
rational numbers in the reduced form $p/q$; this function is discontinuous at every rational
point and continuous at every irrational number.} Bolzano added more excitement in 1834
by giving the first example of a nowhere differentiable continuous function.
As much as works on continuity by Bolzano and Cauchy were ignored, the 
world seem to have started paying attention to the definition of continuity
given in [Darboux, 1875], where its local nature is finally underlined.
\vfill
\eject

{\bf Answers.}

Let's look closer into the questions we posed at the beginning. 
\sju
{\it  A Socratic dialog.}
 
\lqq Do you love your cat?," Arthuro asked.  

\lqq Oh, yes, definitely... and then, no, not at all, most definitely," 
Gwendolin replied. 

\lqq How can that be, aren't you contradicting yourself?" 

\lqq Oh no, no,
I simply do not have a cat!"  exclaimed Gwendolin, and added (philosophically and
pensively):
\lqq I cannot decide on  whether it may be easier not to love a cat 
you do not have than to love a cat you do not have..."

\lqq Hmmm," Arthuro sighed unhappily, for the lack of a definite answer.
\lqq But I will give you a cat, just please give me definitely an answer
whether you love your cat or not."  

\lqq O.K. Arthuro, this would then be a different situation, from the 
previous one... Also, do not forget that
my having a cat may not guarantee that I will be able to decide whether I love
it or not..."

Arthuro is puzzled, but something seems to be happening, for his eyes are 
wide open and he seems to be thinking with intensity...
\sju

When I ask my students whether they love their cats, those who do not have a cat
would never answer the question with a \lqq yes" or \lqq no," but would rather
say that they do not have a cat. This phenomenon of a question with a false 
presumption seems to be well recognized in the folk wisdom; the question of the
same kind is a question of the type \lqq Why do you beat your wife?" 
(while both 
the beating and the wife may not exist). Recently, I  heard a folk child 
riddle: A rooster is standing on top of a roof and it lays an egg.
Which way will the egg go? [The child is supposed to recognize that the
answer is \lqq not applicable," since roosters do not lay eggs.]
\sju

The question of continuity of $1/x$ at $a=0$ is the same as that of
continuity of this function at $a=$ red tomato; or for that matter the
question of Gwendolin loving \lqq her" cat that she does not have. The premise
part of the definition of continuity is not satisfied, since it is vacuous --
there is no function at $0$.
For this reason, saying that $1/x$ is continuous at $0$ is equally (in)valid
as saying that it is continuous at $0$.
The answer to the title question is thus \lqq N/A"
(namely the question of \lqq continuity at 0" is not applicable for the
{\bf function} $1/x$).\footnote{See also [Marde\v si\'c, 1974, p.194]. A more general result holds:
Let $X$ be a Banach algebra with the unity $e$ and let $I$ denote the set of all invertible (regular)
elements in $X$. Then $I$ is an open set and the function $f:I\lra X, f(x)=x^{-1}$ is continuous.  }
In the same manner, we define various kinds of discontinuities of a function (within its domain;
such as done for instance in [Rudin, 1976, p.94]).

As for question number two, we need to find $\lim_{x\to 0} f(x)$. But what 
are the $x$ with $x\to 0$? In fact the only such $x$'s are $x=0$, thus we 
indeed have (4) and (5) satisfied (sufficiently small balls centered at $0$ coincide
with its center). It will always be like this if $a$ is not a limit 
point, if, for instance,  it is an isolated point as in our case.\footnote{See also [Rudin, 1976, p.86]} 
 The function (1)  is continuous (at all points of its domain) and so is $1/x$.
One is tempted to dismiss the \lqq pathological cases" of isolated points 
or non-compact domains in teaching non-mathematicians. (Un)fortunately
discussing continuity only on compact domains would be tantamount to
simply replacing continuity by a stronger notion of uniform continuity.
The field of applied mathematics (aka \lqq the sciences") however largely consists of singularities,
isolated points, non-compact domains and limits that do not exist or
are infinity -- all the staple brushed away in usual calculus courses.

Why does this myth of $1/x$ not being continuous at zero
linger on so persistently with teachers of calculus and authors of 
calculus textbooks? 
If we take aside a small population of people who maintain that incorrect
 mathematics is O.K. (for whatever \lqq higher pedagogical goals"), 
the constraints to open thinking about continuity have historical 
as well as  \lqq objective" roots.

The question of contiguity clearly depends on the space into which the domain 
is embedded. For contiguity purposes, we look at the domain as a part of $\RR$.
If we somehow imagine the domain $(-\infty,0)\cup (0,\infty)$ of $1/x$ not to
be embedded in $\RR$, but standing on its own, we can perhaps imagine that the
graph is contiguous -- that we can draw it in one go, without lifting up the
pencil (after fusing the parts together?)... The same thinking applies to the domain $\{0\}\cup[1,\infty)$ 
of the function (1) in question number two. An interesting point of view is brought
about in [Burgess, 1990], where it is argued that continuity and contiguity
are synonymous in some cases. Thus, if $f:A\lra\RR$ is a function such that
its domain $A=\RR$,  a closed interval,  or is a closed ray in $\RR$, and such that its
graph is closed in $\RR^2$, then $f$ is continuous iff the graph of $f$
is connected in $\RR^2$. 

Another \lqq objective" reason not seeing continuity  for what it is, is the
desire (and sometimes necessity) to extend functions (their domains) so that
the new functions become continuous extensions of the starting functions; this 
is the moment when Arthuro wants to give Gwendolin a cat.
We cannot extend $1/x$ so as to make it continuous at $0$, but we can extend
$x\sin(1/x)$, by defining its value at $0$ to be the limit there, namely $0$.
Notice that if $f:[-1,0)\cup (0,1]\lra \RR, f(x)=1/x$, even
those who say that this $f$ is discontinuous at $0$, would not say that it is
discontinuous at 5, for instance. 

I do however want to pause for a moment from being a devil's advocate: Many
discussions on the nature of the notions of function and continuity did 
spill into the early 20th century and beyond, for instance in the papers of
Borel, Lebesgue, Brower, Baer, etc.; some of these discussions had added
dimension of set-theoretical considerations that led to deep discoveries
in mathematical logic and developments of \lqq new kinds" of mathematics, 
such as intuitionism, etc. It must be said also that there is often a 
considerable delay in adoption and application of fundamental notions
in mathematics. 

\sju

{\bf Conclusions.}

The notion of continuity is subtle, but it had undergone its evolution through
the historical birthing process. The subtleties eluded Cauchy into making a nice
error when he \lqq proved" that a multi-variable function is continuous if it
is continuous in each of its variables.  
The notion of continuity has been 
demystified however, and there is no much reason for dragging the old fog
surrounding it into modern day textbooks and classrooms. 

It is crucial to see functions as ordered triples $f:A\lra B$, for not only
the formula, but both domain and codomain determine what exact properties
the function has. 
This also helps the students understand that  extending a 
function  continuously refers to extension of its domain, etc.
One aesthetic consequence is the quotient rule: if $f$ and $g$ are continuous
then $f/g$ is continuous (the quotient is not a function when the denominator
is zero, thus we do not need to discuss continuity at those points).

While contiguity has its own merits because it is directly related to path
connectedness, it simply is not identical to the notion of continuity. 
As a stronger requirement, it implies continuity. 
We cannot however use the two synonymously -- there are, for instance, 
certain assumptions
about a function that would imply continuity of the function 
(but certainly not contiguity, etc.). 
Once we let go of interpretation
of continuity as contiguity we can concentrate on what continuity actually is:
the commutation of the function with limits as in (4). Why can we plug numbers
into expression whose limits we are finding? Students would invariably plug in
$x=1$ when  finding the limit $\lim_{x\to 1}\root 3\of {x^2+1}$, and an 
opportunity is missed to point out to them that this natural move is
possible because of continuity of the function. In fact plugging in numbers
when finding limits is possible because most of the functions given by formulas
in elementary calculus (e.g. all elementary functions such as algebraic or trigonometric functions) 
 are continuous. Interestingly intuitionists, like Brower, have it that
every function $f:\RR\lra \RR$ is continuous (by their own definitions of functions and
continuity).

Additional discussion on historical development of the concepts  of 
function and continuity
can be found in [Monna, 1972], [Edwards, 1979] and [Hairer and Wanner, 1995].

Acknowledgment: I am indebted to Douglas Klein for pointing out reference 
[Ferraro, 2000] to me and for commenting on  the first version of this paper.
Thanks are also due to Philip Brown for his comments.
\vfill
\eject

{\bf References}
\ju

Jean-B. le Rond d'Alembert: Recherches sur la courbe que forme corde tendu\"e mise
en vibration, {\it M\'emoires de l'Acad\"emie des Sciences de Berlin}, {\bf 3}(1747),
214--219.

B. Bolzano: Rein analytischer Beweis des Lehrsatzes, dass zwischen je zwei
Werthen, die ein entgegengesetztes Resultat gew\"ahren, wenigstens eine reele
Wurzel der Gleichung liege. (Translation: Purely analytical proof of the 
theorem, that between each two
values with the  opposing signs, at least one real root of the
equation lies.) Prag 1817; Ostwald's Klassiker \#153, 1905
[See also O. Stolz: B. Bolzano's Bedeutung in der Geschichte der 
Infinitesimalrechnung. {\it Math. Ann.} {\bf 18}(1879), 255-279].

N. Bourbaki: Les structures fondamentales de l'annalyse, {\it Th\'eorie des 
ensembles} (fascicule de r\'esultats). Paris, 1939.

A.-L. Cauchy: Cours d'analyse alg\'ebrique. \OE uvres s\'erie 2, vol. III., 1821
(Alternatively: Cours D'analyse de l'Ecole Royale Polytechnique, Oeuvres, Ser.2,
Vol. 3, Paris, Gauthier-Villars, 1897).

A.-L. Cauchy: M\'emoire sur les fonctions continues, {\it C. R. Acad. Sci. Paris},
{\bf 18}(1944), 116; see \OE uvres Compl\`etes, Paris, 1882.

G. Darboux: M\'emoire sur les fonctions discontinues, {\it Ann. Ecole Norm.
Sup.}, 2e S\'erie IV (1875), 57--112.

G.L. Dirichlet: \"Uber die Darstellung ganz willk\"urlicher Functionen durch
Sinus- and Cosinusreihen, {\it Repertorium der Physik}, Bd. {\bf I}(1837), 152--174.
In {\it Gesamelte Werke}, Berlin (1889--1897), {\bf I}, 133--160. Also in 
{\it Ostwald's Klassiker}, No. 116, H. Liebmann, Ed., 1900, 3--34.

C.H. Edwards, Jr.: The Historical Development of the Calculus, Springer-Verlag,
New York, 1979.

L. Euler: Opera Omnia, Ser. 2, Vol. 11, Part 2, Series II, 1960.

L. Euler: Introductio in analysin infinitorum, Lausannae: M.M. Bousquet et Soc., 
1748 (or Opera omnia (1), 8--9).

L. Euler: De usu functionum discontinarum in Analysi, {\it Novi Commentarii academiae
scientiarum Petropolitanae}, {\bf 11}(1765), 3-27.

Ferraro, Giovanni: 
Functions, functional relations, and the laws of continuity in Euler, 
{\it Historia Mathematica}, {\bf 27}(2000), 107--132.

J. Fourier: La th\'eorie analytique de la
chaleur, Paris, 1822 (The Analytical Theory of Heat (translated by A. Freeman).
New York, Dover (reprint), 1955).

I. Grattan-Guinness: Bolzano, Cauchy and the `New analysis' of the early
nineteenth century. {\it Arch. Hist. Exact Sci.} {\bf 6}(1969/70), 372--400.

E. Hairer, G. Wanner: Analysis by Its History, Springer-Verlag, New York, 1995.

P.E.B. Jourdain: The origin of Cauchy's conceptions of a definite integral
and of the continuity of a function. {\it Isis} {\bf 1}(1913), 661--703.

S. Marde\v si\'c: Matemati\v cka Analiza. U $n$-dimenzionalnom realnom prostoru.
\v Skolska knjiga, Zagreb, 1974.

A.F. Monna: The concepts of function in the 19th and 20th centuries, in 
particular with regard to the discussions between Baire, Borel and Lebesgue,
{\it Arch. Hist. Exact Sci.}, {\bf 9}(1972), No.1, 57--84.

W. Rudin: Principles of Mathematical Analysis, Third ed., McGraw-Hill, Inc., 
New York, 1976.

D.J. Struik: A Source Book in Mathematics 1200--1800. Harvard University Press,
Cambridge, MA, 1969.

K. Weierstrass: Theorie der analytischen Funktionen, Vorlesung an der Univ. 
Berlin 1874, manuscript (ausgearbeitet von G. Valentin), Math. Bibl. Humboldt
Universit\"at Berlin, 1874.

\end